\documentclass[a4paper,10pt]{article}
\usepackage[hmarginratio=1:1, bottom=3.5cm, top=2.5cm]{geometry}
\usepackage{amssymb,amsthm,amsmath,latexsym,url,graphicx, comment}
\usepackage[english]{babel} 
\usepackage[dvipsnames]{xcolor}
\usepackage[utf8]{inputenc}
\usepackage[sc,osf]{mathpazo}
\usepackage[euler-digits]{eulervm}
\usepackage[T1]{fontenc}
\usepackage{mathtools}
\usepackage{mathtools}
\usepackage{pdfpages}
\usepackage{float}
\usepackage[colorlinks=true,pdftex,unicode=true,linktocpage,bookmarksopen,hypertexnames=false]{hyperref}
\newtheorem{theo}{Theorem}
\newtheorem{prop}[theo]{Proposition}
\newtheorem{lem}[theo]{Lemma}
\newtheorem{cor}[theo]{Corollary}
\newtheorem{rem}[theo]{Remark}

\counterwithout{theo}{section}
\newtheorem{conjecture}[theo]{Question}

\date{}  

\title{On the Sylow Theorem for Skew Braces}
\author{A. Caranti -- I. Del Corso -- M. Di Matteo -- M. Ferrara -- M. Trombetti}
\begin{document}
\maketitle

\begin{abstract} 
\noindent We discuss the (first) Sylow theorem for certain classes of finite skew braces, proving it to hold true when the skew brace is two-sided, bi-skew, right nilpotent, $\lambda$-homomorphic or supersoluble. 
We also show it to hold true for soluble skew braces that are left-nilpotent, and address a number of more specialized settings, proving general Hall-type theorems.
\end{abstract}

\medskip\medskip


\smallskip


\section{Introduction}
A \emph{skew \textnormal(left\textnormal) brace} is a triple $(B, +, \circ)$ where $(B, +)$ and $(B, \circ)$ are (not necessarily abelian) groups that satisfy the {\it skew \textnormal(left\textnormal) distributivity}: $$a \circ (b+c)=a\circ b-a+a \circ c,$$ for all $a, b, c \in B$. We refer to $(B,+)$ as the {\it additive group} and to $(B,\circ)$ as the {\it multiplicative group} of the skew brace $B$. This algebraic structure was introduced in~\cite{GV} as a generalization of the (left) braces defined by Rump in \cite{Rump0}, which correspond to the special case of skew braces whose additive group is abelian. Skew  braces were developed to tackle the challenge of identifying all set-theoretical, non-degenerate solutions to the Yang--Baxter Equation, a fundamental consistency condition that plays a key role in quantum statistical mechanics, underpins the theory of quantum groups, and bridges various disciplines including Hopf algebras, knot theory, braid theory, and Garside theory. 

It is now well established that skew braces sometimes behave like rings (especially when the additive group is abelian) and in other instances like groups, particularly when the additive group is non-abelian. This dichotomy is clearly illustrated in \cite{centnilpo}, where numerous counterexamples to classical group-theoretic results are exhibited; for example, the product of two “nilpotent” ideals need not be “nilpotent”.

Probably, from a group-theorist’s perspective, the most noticeable gap in the current theory is the lack, to date, of an analogue of the Sylow theorems. This is by no means because such results would be irrelevant in this setting: on the contrary, the existence of well-behaved substructures often provides deep insight in many contexts related to the study of solutions of the Yang--Baxter Equation itself (see for example \cite{supersol},\cite{ballester2024soluble},\cite{trombetti2023structure}). The reasons are essentially twofold: there is no a priori reason for such theorems to hold in the setting of skew braces, and, from a computational perspective, identifying counterexamples appears to be difficult (even by using GAP \cite{GAP}). As proof of the difficulty of this task, one may consider that only very recently it has been possible to prove that the only finite skew braces with no proper non-zero sub-skew braces are the ones of prime order (see \cite{ballester2024soluble}). Thus, in a finite simple skew brace, we know by this result that there always exists a sub-skew brace of prime order, but we do not know (for example) if, given a prime dividing the order of the skew brace, there is a sub-skew brace of that order; in other words, we do not know whether an analogue of the Cauchy theorem holds. Despite these discouraging premises, we chose to investigate the analogues of the Sylow theorems for skew braces. Our results show that, in many natural and significant settings, a Sylow-type theorem (and consequently a Cauchy-type theorem) does indeed hold.

\begin{conjecture}
Let $p^n$ be the maximum power of the prime $p$ dividing the order of the finite (soluble) skew brace $(B,+,\circ)$.

Does $B$ contain a sub-skew brace of order $p^n$?
\end{conjecture}

%A preliminary step to prove this conjecture would be proving an analogue of the {\it Cauchy theorem}, that is, the existence, for every prime $p$ dividing the order of the skew brace, of a sub-skew brace of order $p$. Also in this respect, we have some positive results, which makes us hope that at least this result could be true.

\medskip

Actually, our results for (certain classes of) finite soluble skew braces go well beyond proving the existence of~Sy\-low sub-structures. In fact, they provide analogues of what is known as the Hall theorem for finite soluble groups, and that deals with the existence of the so-called {\it Hall subgroups} (analogues of Sylow subgroups for more than one prime). We believe that the following somewhat stronger question might have a positive answer for soluble skew braces.

\begin{conjecture}
Let $\pi$ be any set of primes, and $(B,+,\circ)$ be a finite soluble skew brace. 

Does $B$ contain a sub-skew brace whose order is the maximum $\pi$-number dividing the order of $B$?
\end{conjecture}

(The reason why we often need a strong solubility assumption in our arguments, even for the existence of Sylow subgroups, is explained in Remark \ref{whysoluble}.)

\medskip

The unexplained terminology can be found in Section \ref{preliminaries}. For the reader's convenience, we summarize our main results as follows:
\begin{itemize}
    \item Lemma \ref{lemma:ifkerisone}: this is a key lemma that often allows us to conclude the existence of Sylow and Hall sub-skew braces in favorable contexts. As a consequence of this result, we can almost immediately deduce Sylow and Hall-type theorems for $\lambda$-homomorphic skew braces, right nilpotent skew braces and bi-skew braces (see Corollary \ref{cortot}).

    \item Theorem \ref{leftnilpotent}: a Hall-type theorem for left nilpotent skew braces.

    \item Theorem \ref{bcircbstartb}: a Hall-type theorem for a skew brace $B$ such that both $(B,\circ)$ and $(B\ast B,+)$ have unique (that is, normal) Hall subgroups. As a consequence, we obtain a Cauchy-type theorem when the multiplicative group has a normal Sylow subgroup (see Corollary~\ref{cauchythm}).

    \item Theorem \ref{twosided}: a Hall-type theorem for two-sided skew braces, and so, in particular, for skew braces whose multiplicative group is abelian.

    %\item Theorems \ref{morethanlambda} and \ref{morethanlambda2}: a Cauchy and Hall-type theorem for skew braces $B$ in which~$\lambda_b$ is a skew brace automorphism for every $b\in B$. (But note that these are actually bi-skew braces by Remark \ref{lambaskewbrace}.)

    \item Theorem \ref{thm:supersylow}: a Hall-type theorem for a supersoluble skew brace.

    \item Section \ref{insoluble}: we prove that a well-known family of finite simple skew braces has Hall sub-skew braces.
\end{itemize}

\section{Preliminaries}\label{preliminaries}
A \textit{skew left brace} (resp. {\it skew right brace}) is a triple $(B, +, \circ)$ where $(B, +)$ and $(B, \circ)$ are (not necessarily abelian) groups (respectively, called the {\it additive} and {\it multiplicative} group) such that the {\it skew left distributivity} (resp. {\it skew right distributivity}) holds: $$a \circ (b+c)=a\circ b-a+a \circ c,$$ (resp., $(b+c)\circ a=b\circ a-a+c \circ a$) for all $a, b, c \in B$. Since skew left braces and skew right braces are algebraically equivalent, it is sufficient to restrict our attention to skew {\it left} braces. For simplicity, we will refer to skew left braces simply as skew braces, specifying ‘‘left’’ only when necessary. If $(B,+,\circ)$ satisfies both the skew left and right distributivity, then we say $B$ to be a {\it two-sided} skew brace. If both $(B,+,\circ)$ and $(B,\circ,+)$ are skew braces, then $B$ is said to be a {\it bi-skew brace}.

%The simplest examples of skew braces are the following ones. Given a group $(G,\cdot)$, both $(G,\cdot,\cdot)$ and $(G,\cdot,\cdot^{\operatorname{op}})$ (where $g\cdot^{\operatorname{op}}h=h\cdot g$) are skew braces called respectively the \textit{trivial} and the \textit{almost trivial} skew brace; a trivial skew brace arising from an abelian group is called \textit{abelian}.

Now, let $(B,+,\circ)$ be a skew brace. If $\mathfrak X$ is any class of groups, then $B$ is said to be of $\mathfrak X$ {\it type} if $(B,+)\in\mathfrak X$; thus, for example, $B$ is of {\it abelian type} if $(B,+)$ is an abelian group. In this latter case, we also say that $B$ is a {\it brace} (without the word ‘‘skew’’) because this is the name for the original structures constructed by Rump in \cite{Rump0}.

The map
$$\lambda : a \in (B, \circ) \mapsto \lambda_a \in \operatorname{Aut}(B,+),$$
defined by $\lambda_a(b)=-a+a \circ b$, is a (well-defined) group homomorphism; for every $a,b\in B$, the following identities are verified:
$a+b=a \circ \lambda_a^{-1}(b)$ and $-a=\lambda_a(a^{-1})$,
where $a^{-1}$ denotes the inverse of $a$ in the group $(B, \circ)$. If $\lambda$ is a homomorphism also with respect to the additive group $(B,+)$, then $B$ is said to be {\it $\lambda$-homomorphic} (see \cite{lambdahom}). A third operation, called the {\it star operation}, measures the interaction between the additive and the multiplicative groups of $B$: for each $a,b\in B$, $$a \ast b:=\lambda_a(b) - b=-a + a \circ b-b.$$ The following identities turn out to be useful on many occasions:
$$x \ast (y+z) = x \ast y +y+x \ast z -y\quad\textnormal{and}\quad(x \circ y) \ast z=x \ast (y \ast z)+y \ast z+ x \ast z$$
for all $x, y, z \in B$. Clearly, $a\ast b=0$ for every $a,b\in B$ if and only if $(B,+)=(B,\circ)$; in this case we say that $B$ is {\it trivial} and we can consider $B$ just as a group. In fact, given a group $(G,\cdot)$, we can always realize the trivial skew brace $(G,\cdot,\cdot)$ and the {\it almost trivial} skew brace $(G,\cdot,\cdot^{\operatorname{op}})$, where (where $g\cdot^{\operatorname{op}}h=h\cdot g$). The $\lambda$-action allows to consider the natural semidirect product $(B,+)\rtimes_\lambda(B,\circ)$; a straightforward computation shows that the operation $a\ast b$ corresponds to a commutator of the form
$[(0, a),(b, 0)] = (a \ast b, 0)$ for all $a, b \in B$.

\begin{rem}
{\rm Since every skew brace $(B,+,\circ)$ carries three group structures, we need to distinguish between the usual group-theoretic notations depending on which group operation is involved. This is done by adding a small symbol $+$, $\circ$, or $\lambda$ as a subscript, and by using the adjectives ‘‘additive’’ and ‘‘multiplicative’’ appropriately.  Thus, for example, the additive commutator $a + b - a - b$ is denoted by $[a,b]_+$, while if $a \in (B,+)$ and $b \in (B,\circ)$, we write $[a,b]_\lambda = a - \lambda_b(a)$.}
\end{rem}

It is presumably a well-known fact that bi-skew braces can be characterized as those skew braces for which the $\lambda_a$'s, $a\in B$, are skew braces homomorphisms.

\begin{lem}
Let $B$ be a skew brace. The following are equivalent\textnormal:
\begin{enumerate}
    \item[\textnormal{(1)}] $B$ is a bi-skew brace.
    \item[\textnormal{(2)}] $\lambda_a$ is a skew brace automorphism of $B$ for every $a\in B$.
\end{enumerate}
\end{lem}
\begin{proof}
 %(Here $\lambda_{c}^{-1}$ is the functional inverse, and $c^{-1}$ is the inverse in $(B, \circ)$.)

For all $a, b, c \in B$
    \begin{equation*}
      \lambda_{c}(a \circ b)
      =
      \lambda_{c}(a + \lambda_{a}(b))
      =
      \lambda_{c}(a) + \lambda_{c} \lambda_{a}(b),
    \end{equation*}
    and
    \begin{equation*}
    \lambda_{c}(a) \circ \lambda_{c}(b)
    =
    \lambda_{c}(a) + \lambda_{\lambda_{c}(a)}\lambda_{c}(b).
  \end{equation*}
  Thus the two expressions coincide precisely when for all $a, c \in B$ one has
  \begin{equation*}
    \lambda_{\lambda_{c}(a)}
    =
    \lambda_{c} \lambda_{a} \lambda_{c}^{-1},
  \end{equation*}
  which is, once one translates the notation, condition~4 of
  \cite[Theorem 3.1]{bi}.

  The fact that (1) implies (2) also follows from the fact that if $\lambda$ is the lambda function of the bi-skew brace $(B, +, \circ)$, then the lambda function of $(B, \circ, +)$ is given by $$c\in (B,+) \mapsto \lambda^{-1}_c= \lambda_{c^{-1}}\in \operatorname{Aut}(B,\circ)$$ (see \cite[Subsection 3.1]{bi}).
\end{proof}

\medskip

Now, we briefly review the most important types of sub-structure that one may encounter in a skew brace $(B,+,\circ)$.
\begin{itemize}
    \item \textit{Sub-skew brace}: a subgroup of both $(B,+)$ and $(B,\circ)$. For example, it turns out that $\operatorname{Ker}(\lambda)$ is always a sub-skew brace. As usual in groups and rings, $C\leq B$ means that $C$ is a sub-skew brace of $B$.

    \item ({\it Strong}) \textit{left ideal}: a (normal) subgroup of $(B,+)$ that is also $\lambda$-invariant. Note that a left ideal is always a sub-skew brace, and that every characteristic subgroup of an additive normal subgroup of $B$ is always a strong left ideal.
    
    \item \textit{Ideal}: a strong left ideal of $B$ that is also multiplicative normal subgroup of $B$. If $I$ is an ideal, then $B/I$ inherits the structure of a skew brace with the induced operations on the additive cosets with respect to $I$. As usual in groups and rings, $I\trianglelefteq B$ means that $I$ is an ideal of $B$. The most relevant ideals in a skew brace $B$, (apart from the obvious ones $B$ and $\{0\}$) are probably the following ones:
    \begin{itemize}
        \item The {\it socle} $\operatorname{Soc}(B)$ of $B$: the intersection of $\operatorname{Ker}(\lambda)$ and the centre $Z(B,+)$ of $(B,+)$.
        \item The {\it centre} $Z(B)$ of $B$: the intersection of $\operatorname{Soc}(B)$ and the centre $Z(B,\circ)$ of $(B,\circ)$.
        \item $B\ast B=\langle a\ast b\,:\, a,b\in B\rangle_+$. More in general, if $I$ is any ideal of $B$ and $J$ is any left ideal of $B$, then $I\ast B=\langle i\ast b\,:\, b\in B,\, i\in I\rangle_+\trianglelefteq B$, while $B\ast J=\langle b\ast j\,:\, b\in B,\, j\in J\rangle_+$ is a left ideal of $B$.
    \end{itemize}
    We say that $B$ is {\it simple} if it does not have any non-zero proper ideal.
\end{itemize}

The previous ideals allows us to define natural solubility and nilpotency concepts for a skew brace $(B,+,\circ)$ as follows (see \cite{centnilpo}, \cite{ballester2024soluble}, and \cite{JVV} for more information). The skew brace $B$ is said to be {\it soluble} if it has a finite chain of ideals $$\{0\}=I_0\leq I_1\leq\ldots\leq I_n=B$$ such that $I_{j+1}/I_{j}$ is a trivial skew brace of abelian type for every $0\leq j<n$; in this case, the smallest such non-negative integer $n$ is the {\it soluble length} of $B$. 

Now, we recursively define the {\it upper central series} and the {\it upper socle series}) of $B$ by setting: $Z_0(B)=\{0\}$, \hbox{$Z_1(B)=Z(B)$,} and $Z_{i+1}(B)/Z_i(B)=Z(B/Z_i(B))$ for every $i\geq0$: $\operatorname{Soc}_0(B)$, $\operatorname{Soc}_1(B)=\operatorname{Soc}(B)$, and $\operatorname{Soc}_{i+1}(B)/\operatorname{Soc}_i(B)=\operatorname{Soc}(B/\operatorname{Soc}_i(B))$. Then $B$ is said to be {\it centrally nilpotent} (resp., {\it socle nilpotent}) if $B=Z_n(B)$ (resp., $B=\operatorname{Soc}_n(B)$) for some non-negative integer~$n$; in this case, the smallest such $n$ is the {\it central nilpotency class} (resp., {\it socle nilpotency class}) of $B$. Clearly, every central nilpotent skew brace is socle nilpotent, and every socle nilpotent skew brace is soluble. Also, we observe that when~$B$ is a trivial skew brace the solubility and nilpotency concepts here introduced coincide with the usual ones for groups. For this reason, a soluble skew brace of length $1$ (resp., of length at most $2$) is also said to be {\it abelian} (resp., {\it metabelian}).

Two other nilpotency concept for a skew brace may be defined as follows. In this case, one tries to emulate nilpotency in the context of rings and not groups. We say that $B$ is {\it right nilpotent} if there exists a non-negative integer $n$ such that $$(\ldots(((\underbrace{B\ast B)\ast B\ldots )\ast B}_{\textnormal{$n$ times}}=\{0\};$$ in this case, the smallest such $n$ is the {\it right nilpotency class} of $B$. Note also that the sets appearing in the left-hand side of the above equality are always ideals of $B$. Similarly, $B$ is {\it left nilpotent} if $$\underbrace{B\ast (\ldots \ast (B\ast(B\ast B}_{\textnormal{$n$ times}})))\ldots)=\{0\}$$ for some non-negative integer $n$; in this case, the smallest such $n$ is the {\it left nilpotency class} of $B$. Note that the sets appearing in the left-hand side of the above equality are always left ideals of $B$, but they are not ideals in general. Clearly, every central nilpotent skew brace is both right and left nilpotent, while every socle nilpotent skew brace is right nilpotent. There are many other natural relations between these nilpotency concepts and the nilpotency of the additive/multiplicative groups of the skew brace, but since we do not directly deal with these we refer the interested reader to \cite{centnilpo} and its bibliography. %For our purposes, we only note the following fact.

%\begin{prop}\label{lemtwosidedcentral}
%Let $B$ be a finite two-sided skew brace of order $p^n$ for some prime $p$. Then~$B$ is centrally nilpotent.
%\end{prop}
%\begin{proof}
%It is enough to show that $Z(B)$ is not zero if $B$ is non-zero. Since $G=(B,+)\rtimes_\lambda(B,\circ)$ has prime power order, so $G$ is nilpotent. Set $Z=Z(G)\cap (B,+)$; in particular, $Z$ is a sub-skew brace of $\operatorname{Ker}(\lambda)$. Since the inner multiplicative automorphisms are skew brace automorphisms (see \cite{trap},~Pro\-position 2.3), $Z^g=Z$ for every $g\in B$, so $Z$ is a multiplicative normal subgroup. Consequently, $H=Z\cap Z(B,\circ)\neq\{0\}$. Thus, $H\subseteq Z(B)\neq\{0\}$, and the statement is proved.~\end{proof}

\medskip

Somewhat between nilpotency and solubility, there is the concept of supersoluble skew brace (see \cite{supersol}). The finite skew brace $(B,+,\circ)$ is said to be {\it supersoluble} if it has a finite chain of ideals whose factors have prime order. It turns out for example that every finite skew brace of square-free order is supersoluble (see \cite{supersol}, Corollary 3.9) and that a finite supersoluble skew brace $B$ of order $p_1^{r_1}\ldots p_{n}^{r_n}$, where $p_n<\ldots <p_1$ are primes, has a series of ideals 
$$\{0\}=I_0\leq I_1\leq\ldots \leq I_n=B,$$ where 
$|I_{i+1}/I_i|=p_{i+1}^{r_{i+1}}$ for $0\leq i<n$ (see \cite{supersol}, Theorem 3.11). %In particular, $I_1$ is always a Sylow $p_1$-sub-skew brace of $B$. 
Clearly, every centrally nilpotent skew brace is supersoluble, and every supersoluble skew brace is soluble.

\bigskip

Now, we turn to the main objects of this paper. Let $(B,+,\circ)$ be a finite skew brace, and let~$p$ be a prime. We say that the sub-skew brace~$S_p$ is a {\it Sylow $p$-sub-skew brace} (or also a {\it $p$-Sylow} for short) of $B$ if $|S_p|$ is the largest power of $p$ dividing the order of $B$ (this is equivalent to requiring that $S_p$ is either an additive or a multiplicative Sylow $p$-subgroup of $B$). In a similar fashion, if $\pi$ is a set of primes, we say that the sub-skew brace $R_\pi$ is a {\it Hall $\pi$-sub-skew brace} (or also a {\it $\pi$-Hall} for short) of $B$ if its order is the largest $\pi$-number dividing~$|B|$. Again, this is equivalent to requiring that $R_\pi$ is either an additive or a multiplicative Hall $p$-subgroup. Obviously, if~\hbox{$\pi=\{p\}$,} then a $\pi$-Hall is just a $p$-Sylow. Note also that for a trivial skew brace these concepts coincide with the well-known group-theoretic analogues, and that in the almost trivial case, every additive/multiplicative Hall $\pi$-subgroup is a Hall $\pi$-sub-skew brace. The following result is easily proved and will be used several times without a further notice.

\begin{lem}
Let $(B,+,\circ)$ be a finite skew brace, and let $S$ be a Hall $\pi$-sub-skew brace for some set $\pi$ of primes. If $I$ is any ideal of $B$, then $S+I/I$ and $S\cap I$ are Hall $\pi$-sub-skew braces of $B/I$ and $I$, respectively.
\end{lem}

\noindent Obviously, a finite skew brace will said to be a {\it $\pi$-skew brace} if its order is a $\pi$-number, or equivalently if it coincides with its the unique Hall $\pi$-sub-skew brace.

\smallskip

Clearly, every finite skew brace of nilpotent type has (unique) Hall $\pi$-sub-skew braces for every choice of the set $\pi$ of primes because the additive Hall $\pi$-subgroups are additive characteristic subgroups and so sub-skew braces (see \cite{referee1} for some additional piece of information on the structure of the Hall $\pi$-sub-skew braces in this case). More generally, a finite skew brace has a (unique) Hall $\pi$-sub-skew brace provided that its additive group has a normal Hall $\pi$-subgroup (because it is a characteristic additive subgroup and so a left ideal). Thus, for example, if $p$ is the largest prime dividing the order of a finite supersoluble skew brace $B$, then $B$ has (a unique) Sylow $p$-sub-skew brace.

We also remark that if $(B,+,\circ)$ be a finite skew brace of order a power of a prime $p$, then it is easy to see that $B$ has a sub-skew brace of order $p$ (see also the main theorem of \cite{ballester2024soluble}). Therefore, if we can prove the existence of a (non-zero) Sylow $p$-sub-skew brace in a certain environment, then we can automatically also prove the Cauchy theorem for the prime $p$. %With respect to the Cauchy theorem, we would like to note that one day before we submitted this manuscript, the manuscript "On a Cauchy theorem for finite skew braces" appeared on ArXiv, in which the main results are the Cauchy theorem for bi-skew and two-sided skew braces.

\medskip

One of the tools we will employ to prove our main result is the opposite skew brace. Recall that for any skew brace $(B,+,\circ)$, it is possible to define another skew brace which goes under the name of {\it opposite} skew brace of $B$, and is denoted by $B^{\operatorname{opp}}$. This skew brace is obtained by replacing the additive group $(B,+)$ with its {\it opposite} group $(B,+^{\operatorname{opp}})$, which is defined by: $a+^{\operatorname{opp}}b=b+a$ for every $a,b\in B$. The concept of opposite skew brace has been introduced in \cite{opposite} (but see also \cite{campedel}), and since then it has proved to be very useful in the study of skew braces and their applications. From our point of view, the most relevant feature of the opposite skew brace is that it preserves sub-skew braces and (strong left) ideals. The following statement is certainly well known, but we were unable to find it in the literature. For completeness, we include a short proof.

\begin{lem}\label{lemma5}
Let $(B,+,\circ)$ be a skew brace, and let $C$ be a sub-skew brace \textnormal(resp., strong left ideal, or  ideal\textnormal) of $B$. 
Then $C$ is a sub-skew brace \textnormal(resp., strong left ideal, or ideal\textnormal) of $B^{\operatorname{opp}}$.
\end{lem}
\begin{proof}
The fact that $C$ is a sub-skew brace of $B^{\operatorname{opp}}$ is obvious because the (normal) subgroups of $(B,+)$ are also (normal) subgroups of $(B,+^{\operatorname{opp}})$. Suppose $C$ is a (strong left) ideal of $B$. We only need to prove that $C$ is $\lambda$-invariant. Thus, let $a\in B$ and $b\in C$. We have that $$\lambda_a(b)^{\operatorname{opp}}=-a+^{\operatorname{opp}}a\circ b=a\circ b-a=a-a+a\circ b-a=a+\lambda_a(b)-a$$ belongs to $C$ because $C$ is a strong left ideal of $B$.
\end{proof}

\section{Sylow and Hall Theorems}

We start this section with a key lemma that often allows us to conclude the existence of Sylow and Hall sub-skew braces whenever the kernel of $\lambda$ behaves.

\begin{lem}\label{lemma:ifkerisone}
Let $p$ be a prime \textnormal(resp. $\pi$ be a set of primes\textnormal), and let $(B,+,\circ)$ be a finite  skew brace \textnormal(resp., a finite skew brace with $(B,+)$ and $(B,\circ)$ soluble\textnormal) such that $(B,\circ)/\operatorname{Ker}(\lambda)$ or~$(B,\circ)/\operatorname{Ker}(\lambda^{\operatorname{opp}})$ is a $p$-group \textnormal(resp. a $\pi$-group\textnormal). If $P$ is any left ideal of $B$ of $p$-power order \textnormal(resp. of order a $\pi$-number\textnormal), then~$B$ has a Sylow $p$-sub-skew brace $P_2$ \textnormal(resp. a~Hall $\pi$-sub-skew brace\textnormal) containing $P$. \textnormal(Actually, $P_2$ turns out to be a left ideal of $B$ or of $B^{\operatorname{opp}}$.\textnormal)
\end{lem}
\begin{proof}
We only prove the case for a prime $p$, the proof for the case of a set of primes $\pi$ being totally analogous (the only difference is the use of the theorems of Hall instead of those of~Sy\-low). Since the sub-skew braces of $B$ and its opposite are the same by Lemma~\ref{lemma5}, it is enough to prove the case when $(B,\circ)/\operatorname{Ker}(\lambda)$ is a $p$-group. Let $G=(B,+)\rtimes_\lambda(B,\circ)$. Now, since $P$ is a~\hbox{$p$-sub-skew} brace of~$B$, $S=(P,+)\rtimes_\lambda (P,\circ)$ is a~\hbox{$p$-sub}\-group of $G$. By the hypotheses,~$S$ can be extended to a~\hbox{$p$-sub}\-group of the form $S_1=(P,+)\rtimes_\lambda (P_1,\circ)$, where $(P_1,\circ)$ is a~Sy\-low~\hbox{$p$-sub}\-group of $(B,\circ)$ containing~$(P,\circ)$, and subsequently to a~Sy\-low~\hbox{$p$-sub}\-group of $G$, which is necessarily of the form $S_2=(P_2,+)\rtimes_\lambda (P_1,\circ)$, where $(P_2,+)$ contains~$(P,+)$. Consequent\-ly,~\hbox{$\lambda_a(P_2)=P_2$} for eve\-ry~\hbox{$a\in P_1$.} But $(B,\circ)=P_1\circ\operatorname{Ker}(\lambda)$, and so~$P_2$ is~\hbox{$\lambda$-in}\-va\-riant. Therefore,~$P_2$ is the required~Sy\-low~\hbox{$p$-sub-skew brace.}
\end{proof}

\begin{cor}\label{cortot}
Let $p$ be a prime \textnormal(resp. $\pi$ be a set of primes\textnormal), and let $(B,+,\circ)$ be a finite  skew brace \textnormal(resp., a finite skew brace with $(B,+)$ and $(B,\circ)$ soluble\textnormal). Suppose that one of the following conditions holds\textnormal:
\begin{itemize}
    \item[\textnormal{(1)}] $B$ is a bi-skew brace.\quad \textnormal{(2)} $B$ is right nilpotent. 
    \quad \textnormal{(3)} $B$ is $\lambda$-homomorphic.
    %\item[\textnormal{(2)}] $B$ is $\lambda$-homomorphic.
    %\item[\textnormal{(3)}] $B$ is right nilpotent.
\end{itemize}
\noindent If $S$ is any left ideal of $B$ whose order is a $p$-number \textnormal(resp., a $\pi$-number\textnormal), then $B$ has a Sylow $p$-sub-skew brace \textnormal(resp. a Hall $\pi$-sub-skew brace\textnormal) containing $S$.
\end{cor}
\begin{proof}
(1) \quad In this case $\lambda$ is an antihomomorphism of $(B,+)$ (see \cite{bi}, Theorem 3.1 (3)), so that the kernel $K = \operatorname{Ker}(\lambda)$ is an ideal of $B$. If $K = \{ 0 \}$, then $B$ is an almost trivial skew brace (see \cite{bi}, The\-o\-rem~3.1~(7)), for which the results hold true trivially. Let $B$ be a minimal counterexample, so that $K \ne \{ 0 \}$. By the minimality of $B$ we have that $B/K$ is a $\pi$-skew brace, and thus the results follow from~Lem\-ma~\ref{lemma:ifkerisone}.

\medskip

\noindent(2)\quad We use induction on the right nilpotency class $n$ of $B$. Clearly, if $n\leq 1$, then the statement is trivial. Suppose $n>1$, and let $I$ be an ideal of $B$ such that $I\ast B=\{0\}$ and $B/I$ is right nilpotent of right nilpotency class $n-1$. By induction $B/I$ has a~Hall~\hbox{$\pi$-sub-skew} brace, so without loss of generality we may assume $|B/I|$ is a $\pi$-number. Since $I\subseteq\operatorname{Ker}(\lambda)$, we may apply Lem\-ma~\ref{lemma:ifkerisone} to conclude the proof.

\medskip

\noindent(3)\quad This follows from~\cite{stefanellotrap}, Theorem~3.13, and (2).\end{proof}

\medskip

Left nilpotent skew braces cannot be dealt with as the right nilpotent ones, the reason being that the sub-structures that define them are usually only left ideals and not ideals. 

\begin{rem}\label{whysoluble}
{\rm The proof of Theorem \ref{leftnilpotent} will also make it clear why (a strong) solubility assumption for skew braces is necessary for our arguments to work. In fact, we need a minimal ideal $I$ of the skew brace $B$ to be abelian, a fact that is only possible if $B$ is  soluble in the sense of \cite{ballester2024soluble}. In this circumstances for example every characteristic subgroup of {\it the group} $I$ is an ideal of $B$, being a characteristic subgroup of both groups, so $I$ must also be characteristically simple and consequently a $q$-group for some prime $q$. These are standard arguments that we will often employ.}
\end{rem}

\begin{lem}\label{lemincredible}
Let $(B,+,\circ)$ be a finite, left nilpotent skew brace.  Let $(B,+)=(I,+)+(J,+)$, where~$I$ and~$J$ are strong left ideals of $B$ of coprime orders, and $I\trianglelefteq B$. Then $J\trianglelefteq B$.
\end{lem}
\begin{proof}
We only need to prove that $J$ is multiplicatively normal. Let $B$ be a minimal counterexample to the statement. Since $B$ is left nilpotent, $I+(J\ast J)$ is a proper ideal of $B$, so $J\ast J$ is an ideal of $I+(J\ast J)$, and hence also an ideal of $B$. By minimality of $B$, $J\ast J=\{0\}$. Now, since~$B$ is left nilpotent, $(I,+)\rtimes_\lambda(J,\circ)$ is a group in which $(I,+)$ and $(J,\circ)$ are coprime, and $(J,\circ)$ stabilizes a series of $(I,+)$. This implies that $J\ast I=\{0\}$, so $J\ast B=\{0\}$. Thus, $J\leq \operatorname{Ker}(\lambda)$ and hence $[J,I]_\circ=\{0\}$, which means that $J$ is multiplicatively normal, a contradiction.
\end{proof}

\begin{theo}\label{leftnilpotent}
Let $\pi$ be a set of primes, and let $(B,+,\circ)$ be a finite soluble skew brace which is left nilpotent. Then the Hall $\pi$-sub-skew braces of $B$ coincide with the Hall $\pi$-subgroups of $(B,\circ)$.
\end{theo}
\begin{proof}
Let $B$ be a minimal counterexample to the statement, and let $I$ be a minimal ideal of~$B$ contained in $B\ast B$. By minimality of~$B$, we have that $|B/I|$ is a $\pi$-number. Moreover, since $B$ is soluble, $I$ is an abelian $q$-group for some prime $q\not\in\pi$. Let $P$ be any Hall $\pi$-subgroup of~$(B,\circ)$, and extend it to a Hall $\pi$-subgroup $P_1\rtimes_\lambda P$ of $G=(B,+)\rtimes_\lambda(B,\circ)$. By hypothesis, $B$ is left nilpotent, so $P$ acts on $I$ with respect to $\lambda$ as a nilpotent group of automorphisms (that is, it stabilizes a finite series of $(I,+)$), and hence $P\ast I=\{0\}$.

Let $a\in I$, $b\in P_1$, and $c\in P$. Then $b+a-b=a_1$ for some $a_1\in I$. Since $I\leq B\ast B$, $P\ast I=\{0\}$, so conjugating the previous equality by $c$ in $G$ yields \hbox{$\lambda_c(b)+a-\lambda_c(b)=a_1$,} so~\hbox{$-b+\lambda_c(b)=[-b,c]_\lambda$} additively centralizes $a$. The arbitrariness of $a$, $b$ and $c$ shows that~\hbox{$P\ast P_1$} additively centralizes $I$.

It follows that $(B\ast B,+)$ is the direct sum of $(I,+)$ and $(P\ast P_1,+)$. By Lemma \ref{lemincredible}, $P\ast P_1$ is an ideal of $B\ast B$, and so also an ideal of $B$. Thus, the minimality of $B$ yields that $P\ast P_1=\{0\}$. But then $P\ast B=P\ast (P_1+I)=\{0\}$, so $P\ast P=\{0\}$, and hence~$P$ is a Hall $\pi$-sub-skew brace of~$B$.
\end{proof}

\medskip

We remarked in Section \ref{preliminaries} that the Hall-type theorems trivially hold when the additive group has a unique (and so normal) Hall subgroup. Our next result shows that such hypothesis on the additive group can be weakened if we require a similar one on the multiplicative group.

\begin{theo}\label{bcircbstartb}
Let $\pi$ be a set of primes and $(B,+,\circ)$ a finite skew brace. If both $(B,\circ)$ and $(B\ast B,+)$ have normal Hall $\pi$-subgroups $P$ and $P_1$, respectively, then $P$ is a Hall $\pi$-sub-skew brace of~$B$.
\end{theo}
\begin{proof}
Suppose the statement is false, and let $B$ be a minimal counterexample to the statement; in particular, $B\ast B\neq\{0\}$ otherwise $B$ would be  a group and the statement trivially holds. Let $I$ be a minimal ideal of $B$ contained in $B\ast B$. Since $B$ is soluble, $I$ is a trivial skew brace of abelian type, and so its characteristic subgroups are ideals of $B$. This implies that~$(I,+)$ is an elementary abelian~\hbox{$q$-group} for some prime $q$. By minimality of $B$, we have that $|B/I|$ is a~\hbox{$\pi$-num}\-ber and (obviously) that $q\not\in\pi$.
By hypothesis, $(B\ast B,+)$ contains a characteristic Hall~\hbox{$\pi$-sub}\-group~$P_1$. But then $P_1$ is a sub-skew brace (and in particular a subgroup of $(B,\circ)$), which, by order considerations, coincides with the only Hall $\pi$-subgroup of $(B\ast B,\circ)$. Thus,~$P_1$ is an ideal of $B$. Now, $B\ast B$ is the direct product (as a skew brace) of $I$ and $P_1$, so if we factor~$P_1$ out, then we may assume that $I=B\ast B$.

Let $P$ be the unique Hall $\pi$-subgroup of $(B,\circ)$. In the natural semi-direct product $G=(B,+)\rtimes_\lambda(B,\circ)$, extend $P$ to a Hall $\pi$-subgroup $R$ of $G$. Since $(B,+)$ is normal in $G$, there is a Hall~\hbox{$\pi$-sub}\-group $P_2$ of $(B,+)$ such that $R=(P_2,+)\rtimes_\lambda (P,\circ)$ (use Schur Zassenhaus). Thus, if $a\in P_2$ and $x\in P$, then $x\ast a\in P_2$; but $B/I$ is a trivial skew brace, so $x\ast a\in I$, and hence $x\ast a\in P_2\cap I=\{0\}$. Therefore,~\hbox{$P\ast P_2=\{0\}$.}

If $b\in (B,\circ)$, then the action by conjugation in $G$ on $R$ gives $$^bR=\lambda_b(P_2)\rtimes_\lambda P,$$ which shows that $P\ast\lambda_b(P_2)=\{0\}$. Let $U$ be the subgroup of $(B,+)$ generated by the $\lambda_b(P_2)$'s with $b\in B$. Then $P\ast U=\{0\}$. By its very definition,~$U$ (and consequently $U\cap I$) is $\lambda$-invariant, so $U$ is a sub-skew brace. Since $P_2$ is an additive subgroup of $U$, the order of $U$ is divided by~\hbox{$|P_2|=|P|$}, and this means that the multiplicative Hall $\pi$-subgroup of $U$ must coincide with~$P$ (being the latter the unique Hall $\pi$-subgroup of $(B,\circ)$). Consequently, $P\ast P=\{0\}$ and so $P$ is a~Hall~\hbox{$\pi$-sub-skew} brace of $B$.%, the final contradiction.
%We claim that $U=B$. Clearly, $U\cap I$ is centralized by $(I,+)$. But since $U$ contains $P_2$, which is an additive Sylow $p$-subgroup of $(B,+)$, and also a complement to $(I,+)$, we have that $(U\cap I,+)$ is normalized by $(P_2,+)$ and so even by $(B,+)=(P_2+I,+)$. 
\end{proof}

\begin{cor}\label{cauchythm}
Let $p$ be a prime dividing the order of the skew brace $(B,+,\circ)$. If $P$ is a multiplicative normal~Sy\-low~\hbox{$p$-sub}\-group of $B$, then $B$ has a sub-skew brace of order $p$.
\end{cor}
\begin{proof}
Suppose the statement is false, and let $B$ be a minimal counterexample to the statement. Let $I$ be a minimal ideal of $B$, so $(I,+)$ is an elementary abelian~\hbox{$q$-group} for some prime $q\neq p$. By minimality, $|B/I|=p$, so $B\ast B\leq I$. The contradiction now follows from Theorem \ref{bcircbstartb}.
\end{proof}

\medskip

Next, we deal with two-sided skew braces, and here the argument makes great use of some non-trivial facts about the their automorphism group.

\begin{theo}\label{twosided}
Let $p$ be a prime \textnormal(resp. $\pi$ be a set of primes\textnormal), and let $(B,+,\circ)$ be a finite  \textnormal(resp., finite soluble\textnormal) two-sided skew brace. Then the~Sylow $p$-sub-skew braces \textnormal(resp., Hall $\pi$-sub-skew braces\textnormal) of $B$ coincide with the  Sylow $p$-subgroups \textnormal(resp., Hall $\pi$-subgroups\textnormal) of $(B,\circ)$.
%If $S$ is a~\hbox{$\pi$-sub}\-group of $(B,\circ)$, then there exists a Hall $\pi$-sub-skew brace of $B$ containing $S$. In particular,~$B$ has~Hall $\pi$-sub-skew braces, and the Hall $\pi$-sub-skew braces of $B$ are multiplicatively conjugated.
\end{theo}
\begin{proof}
We prove the result for a set $\pi$ of primes, the case of a single prime $p$ being easier. Since the multiplicative inner automorphisms of $B$ are skew braces automorphisms (see \cite{trap},~Pro\-position 2.3) and $(B,+)$ is soluble, we only need to prove the existence of the Hall~\hbox{$\pi$-sub-skew} braces (because of the well-known theorems of Philip Hall). Thus, let $B$ be a minimal counterexample to the existence, and set $H=(B\ast B)\cap (B\ast^{\operatorname{opp}} B)$. Then $H$ is an ideal of $B$ and~$(H,+)$ is abelian by Theorem 4.3 of \cite{trap}. 

%By Lemma 4.25 and Theorem 4.26 of \cite{trap}, we have that $H$ is centrally nilpotent.

Suppose first that $H=\{0\}$. Clearly, neither $B\ast B$ nor $B\ast^{\operatorname{opp}}B$ can be $\{0\}$ (because otherwise~$B$ has Hall $\pi$-sub-skew braces). By minimality of $B$, both the orders of $B/B\ast B$ and  $B/B\ast^{\operatorname{opp}}B$ are~\hbox{$\pi$-numbers.} Thus, $|B|$ is a~\hbox{$\pi$-number} itself, and we have a contradiction.

Thus, $H\neq\{0\}$. Let $I$ be a minimal ideal of $B$ contained in $H$. Then~$(I,+)$ is an elementary abelian $q$-group for some prime $q$. By minimality, the statement holds for $B/I$, so we have that the order of $B/I$ is a $\pi$-number and $q\not\in\pi$. Also $(H,+)$ is abelian, so it contains a unique~Hall~\hbox{$\pi$-sub}\-group~$R$, which is then an ideal of $H$ (see \cite{trap}, Corollary 2.4). But this means that $H=I\times R$, so $R$ is actually both additively and multiplicatively characteristic, implying that $R$ is an ideal of $B$. By minimality of $B$, $R=\{0\}$. Thus, $H=I$ is contained in $Z(B\ast B,+)$ by~\cite{trap},~The\-o\-rem~4.3. Let $J/I$ be a minimal ideal of $B/I$ contained in $B\ast B/I$, if any. Then $J$ is an ideal of $B$ of nilpotent type. Again, if $S$ is a~Hall~\hbox{$\pi$-sub}\-group of $(J,+)$, then $S$ is an ideal of~$B$ which must then be $\{0\}$ by minimality of~$B$. Therefore~\hbox{$I=H=B\ast B$.}

%Assume that the centre $Z$ of $(B,+)$ is non-zero. Then the Sylow subgroups of $Z$ are ideals of~$B$, being additive characteristic subgroups (see \cite{trap}, Corollary 2.4). Consequently,~$Z$ must be a~\hbox{$q$-group,} and hence $Z\leq I$. By minimality of $B$, we have that $Z=I$, so $(B,+)$ is nilpotent. But then the Sylow $p$-subgroup of $(B,+)$ is an ideal of $B$, and we have a contradiction. Thus, the centre of $(B,+)$ is trivial.

%Let $C=(B,+_{\textnormal{opp}},\circ)$ be the opposite skew brace of $(B,+,\circ)$. Then again $I$ is a minimal ideal of $C$, and both $C/I$ and $I$ are trivial braces.

Let $P$ be a Hall $\pi$-subgroup of $(B,\circ)$. Then there is a Hall $\pi$-subgroup $P_1$ of $(B,+)$ such that $P_1\rtimes_\lambda P$ is a Hall $\pi$-subgroup of $G=(B,+)\rtimes_\lambda(B,\circ)$. Clearly, $P\ast P_1\subseteq I\cap P_1=\{0\}$. %Moreover, if $P_1\subseteq P$, then $P_1=P$ (by order considerations) is a Sylow $p$-sub-skew brace of $B$, a contradiction. Thus, $P_1\setminus P\neq\emptyset$.
%-a+_opp ab=ab-a=+a-a+ab-a
Now, let $U$ be the largest subgroup of $(B,+)$ such that $P\ast U=\{0\}$. Clearly, $P_1\leq U$. We claim that $P_1<U$. Suppose not, so $P_1=U$. Now, since $B$ is two-sided, so $a\ast y=0$ implies $a\ast y^{-1}=0$ for every $a,y\in B$. Thus, $P_1^{-1}\subseteq U=P_1$. Similarly, if $a\ast b=a\ast c=0$ for some $a,b,c\in B$, then $a\ast (b\circ c)=0$. Consequently, $P_1\circ P_1\subseteq U=P_1$. This implies that $P_1$ is a multiplicative subgroup and so a Hall $\pi$-sub-skew brace of $B$, a contradiction.

This last contradiction shows that $U>P_1$. Put $V=U\cap I$. Then $V=U\cap I$ is not $\{0\}$ and is additively normal because $U$ contains an additive complement to $I$. Now, $P\ast V=\{0\}$, and so also $B\ast V=\{0\}$ because $(B,\circ)=I\circ P$ and $I\ast V=\{0\}$. 

Since the multiplicative inner automorphisms are skew brace automorphisms, we have that $B\ast \big(x\circ V\circ x^{-1}\big)=\{0\}$ for all $x\in (B,\circ)$, so $x\circ V\circ x^{-1}\subseteq U$ and hence $V=x\circ V\circ x^{-1}$. Therefore~$V$ is multiplicatively normal. Since $B\ast V=\{0\}$, so $V$ is an ideal of $B$, which must coincide with the minimal ideal $I$. Thence $B\ast I=\{0\}$.

Finally, if $a,b$ are multiplicative $\pi$-elements of $B$, then $\lambda_a(b)=b+u$ for some $u\in I$ (because $I=B\ast B$). But \hbox{$\lambda_a(u)=u$} and $u$ has additive order dividing $q$, so $u=0$ (because $q\not\in\pi$), and consequently \hbox{$a\ast b=0$.} This shows that every Hall $\pi$-subgroup of $(B,\circ)$ is a~(Sy\-low~$\pi$-)sub-skew brace of $B$ and gives the last contradiction.
\end{proof}

\begin{cor}
Let $\pi$ be a set of primes and let $(B,+,\circ)$ be a finite soluble skew brace such that $(B,+)$ is abelian. Then the~Hall $\pi$-sub-skew braces of $B$ coincide with the  Hall $\pi$-subgroups of $(B,\circ)$.
%If $S$ is a~\hbox{$\pi$-sub}\-group of $(B,\circ)$, then there exists a Hall $\pi$-sub-skew brace of $B$ containing $S$. In particular,~$B$ has~Hall $\pi$-sub-skew braces, and the Hall $\pi$-sub-skew braces of $B$ are multiplicatively conjugated.
\end{cor}

The argument we used for two-sided skew braces can be somewhat modified to prove a Hall-type theorem whenever the $\lambda_b$'s are skew brace automorphisms instead of only automorphisms of the additive group.

\begin{theo}\label{morethanlambda}
Let $\pi$ be a set of primes and $(B,+,\circ)$ a finite soluble skew brace such that $\lambda_b$ is a skew brace isomorphism for every $b\in B$. If $Q$ is any left ideal of $B$, then $Q$ is contained in a Hall $\pi$-sub-skew brace of $B$.
\end{theo}
\begin{proof}
Let $B$ be a minimal counterexample to the statement, and let $I$ be a minimal ideal of $B$. Then~$(I,+)$ is an elementary abelian $q$-group for some prime $q$. By minimality of $B$, we have that $B/I$ has order $p^n$ for some $p\neq q$. 

Let $P$ be a Sylow $p$-subgroup of $(B,\circ)$. Then there is a Sylow $p$-subgroup $P_1$ of $(B,+)$ such that $P_1\rtimes_\lambda P$ is a Sylow $p$-subgroup of $G=(B,+)\rtimes_\lambda(B,\circ)$. Thus, $\lambda_x(P_1)=P_1$ for every $x\in P$. Since $\lambda_x$ is a skew brace isomorphism for every $x\in P$, we have that $\lambda_x(P_1)=P_1$ for every $x\in \lambda_b(P)$ with $b\in P$. 

Let $U$ be the normalizer of $(P_1,+)$ in $(B,\circ)$. Clearly, $U$ contains every $\lambda_x(U)$ with $x\in P$ because $\lambda_x$ is a skew brace isomorphism, and we may write $U=P\circ (I\cap U)$. Moreover, since~$P$ is not a sub-skew brace, $I\cap U\neq\{0\}$. If $c\in I\cap U$, then $\lambda_c(P_1)=P_1$. On the other hand, $(I,+)\rtimes_\lambda(I,\circ)$ is a normal Sylow $q$-subgroup of $G$, so $\lambda_c(u)=u$ for every $u\in P_1$. Therefore $c\ast P_1=\{0\}$. But obviously, $c\ast I=\{0\}$, and hence $c\in K=\operatorname{Ker}(\lambda)$.  It follows that $I\cap U\leq K$, which means that $I\cap U=K\cap I$ because $K\leq U$. %Since $B$ is a bi-skew brace, so $K$ is a left-ideal of $B$, and hence $U=P+(K\cap I)$

Now, if $a\in P$, $b\in I\cap U\leq K$ and $u\in U$, then $\lambda_{a\circ b}(u)=\lambda_a(u)\in U$. Thus, $\lambda_u(v)\in U$ for every $u,v\in U$, and  consequently, $u+v=u\circ\lambda_{u^{-1}}(v)\in U$ for every $u,v\in U$. By choosing $w=\lambda_u(u^{-1})$, we see that $u+w=0$. Therefore, $U$ is a sub-skew brace of $B$. Since $|P|$ divides the order of $|U|$,  we have by minimality that $U=B$, which implies that
 $I\subseteq K$. Lemma \ref{lemma:ifkerisone} provides the final contradiction. 
\end{proof}

\medskip

%A proof similar to that of Theorem \ref{morethanlambda} allows us to prove the following result.

%\begin{theo}\label{morethanlambda2}
%Let $p$ be a prime dividing the order of a finite skew brace $(B,+,\circ)$, and assume that $\lambda_b$ is a skew brace isomorphism for every $b\in B$. Then $B$ has a sub-skew brace of order $p$.
%\end{theo}

%\begin{rem}\label{lambaskewbrace}
%{\rm In fact, {\it skew braces $B$ whose $\lambda$-maps are skew brace isomorphisms are bi-skew braces}, so the previous results could also be deduced from Corollary \ref{cortot}. However, we have chosen to include the proofs, as they provide a different perspective from the bi-skew brace case.}
%\end{rem}
%\begin{proof}
%Let $a,b,c \in B$. Then one has that
%\[
%\lambda_a(b) + \lambda_a\lambda_b(c)
%= \lambda_a(b \circ c)
%= \lambda_a(b) \circ \lambda_a(c)
%= \lambda_a(b) + \lambda_{\lambda_a(b)}\lambda_a(c).
%\]
%This shows that
%\[
%\lambda_a \lambda_b = \lambda_{\lambda_a(b)} \lambda_a.
%\]
%However, we know that
%\[
%\lambda_a \lambda_b = \lambda_{\lambda_a(b)}\lambda_{\rho_b(a)}.
%\]
%Combining both identities gives us that
%\[
%\lambda_a = \lambda_{\rho_b(a)}.
%\]
%Filling in the formula for $\rho$ delivers us that
%\[
%\lambda_a = \lambda_{(a^{-1}+b)^{-1} \circ b}.
%\]
%Hence,
%\[
%\lambda_a \lambda_{b^{-1}} = \lambda_{(a^{-1}+b)^{-1}},
%\]
%or equivalently that
%\[
%\lambda_{b \circ a^{-1}} = \lambda_{a^{-1}+b}.
%\]

%This shows that $(B,+,\circ)$ is a $\lambda$-anti-homomorphic skew brace. Thus, by Proposition~3.7 of~\cite{referee2}, $(B,+,\circ)$ is a bi-skew brace.
%\end{proof}

In order to deal with supersoluble skew braces, the following results are crucial. The first provides a complete description of automorphisms of the semidirect product of two groups under specific conditions, while the second is a special case of~\cite[Proposition
  2.24]{campedel}, which for reader's convenience we restate for skew {\it left} braces. %Recall  that if $(G,\cdot)$ is a group, then a {\it lambda function} is....
%\begin{theo}[\protect{Special case of \cite[Theorem 1]{Curran}}]
%  \label{thm:semidirect}
%
%  Let $G = Q \ltimes M$ be a semidirect product, with $M$ abelian, and 
%  characteristic in $G$.
%  
%  Then $\operatorname{Aut}(G)$ can be described in a natural way via the set of
%  matrices
%  \begin{equation}
%    \label{eq:autgrp}
%    \begin{aligned}
%      {\Bigg\{}\ 
%      \begin{bmatrix}
%        d & b\\
%        0 & a
%      \end{bmatrix}
%      :\ &
%      a \in \operatorname{Aut}(M), d \in \operatorname{Aut}(Q),
%      \\&
%      \text{where for $h \in M$ and $k \in Q$ we have}
%      \\&
%      (h^{k})^{a} = (h^{a})^{k^{d}} \text{, for $h \in M, k \in Q$},
%      \\&
%      b : Q \to M,
%      %\\&
%      (x y)^{b} = (x^{b})^{y^{d}} (y^{b})
%      \text{, for $x, y \in Q$}
%      \ {\Bigg\}}.
%    \end{aligned}
%  \end{equation}
%\end{theo}

%In other words, an automorphism
%\begin{align*}
%  \theta
%  =
%  \begin{bmatrix}
%    d & b\\
%    0 & a
%  \end{bmatrix}
%\end{align*}
%acts as
%\begin{align*}
%  (k h)^{\theta}
%  =
%  k^{d} k^{b} h^{a},
%\end{align*}
%for $k \in Q$ and $h \in M$.\\

\begin{theo}[\protect{Special case of \cite[Theorem 1]{Curran}}]
\label{thm:semidirect}

Let $G = M\rtimes P$ be a semidirect product, where $M$ is abelian and characteristic
in $G$. %The action of $Q$ on $M$ is given by conjugation ${}^{k}h = k h k^{-1}$.
Then $\operatorname{Aut}(G)$ can be described in a natural way via the set of
matrices
\[
\Bigg\{
\begin{bmatrix}
a & b\\
0 & d
\end{bmatrix}
:\ 
\begin{array}{l}
a \in \operatorname{Aut}(M),\ d \in \operatorname{Aut}(P),\\[2pt]
a({}^{k}h) = {}^{d(k)}a(h)
\quad \text{for all } h \in M,\ k \in P,\\[4pt]
b : P \to M \text{ satisfies } 
b(xy) = b(x)\,{}^{d(x)}b(y)
\quad \text{for all } x,y \in P
\end{array}
\Bigg\}.
\] Here, an automorphism
\[
\theta=
\begin{pmatrix}
a & b\\
0 & d
\end{pmatrix}
\]
acts as $\theta(hk)= a(h)\, b(k)\, d(k)$, for
$h\in M$, and $k\in P$.
\end{theo}

%\begin{lem}
%  \label{lemma:ifkerisone}
%  Let $(B,\cdot,\circ)$ be a finite right skew brace, such that $(B,\cdot)=QM$, where $Q$ and $M$ are respectively a $q$-Sylow subgroup and a $p$-Sylow subgroup and $M$ is a minimal ideal in $B$ of order $p$.
%  If $\ker(\gamma) \ge M$ or $\ker(\overline{\gamma}) \ge M$, then $B$ has a $q$-Sylow sub-skew brace.
%\end{lem}
%\begin{proof}
%  Since the sub-skew braces of $B$ and its opposite are the same by
%  Lemma~\ref{lemma:opp}\eqref{item:same}, it is 
%  enough to prove the case when $\ker(\gamma) \ge M$.
%
%  The group $\gamma(B)$ of automorphisms of $(B,\cdot)$ has $q$-power
%  order. Therefore in its action on the set 
%  of the Sylow $q$-subgroups of $(B,\cdot)$, which has cardinality
%  $\null\equiv 1 \pmod{q}$, it will fix a
%  Sylow $q$-subgroup, say $Q$.
%
%  In particular, $Q$ is thus invariant
%  under $\gamma(B)$, and is thus the Sylow $q$-sub-skew brace we are
%  looking for.
%\end{proof}

%\begin{rem}
%  As a comment  the last statement in the proof, note that since
%  $(B,\cdot) = Q M$, so that for $x \in M$ and $y \in Q$ we have 
%  \begin{align*}
%    \gamma(x y)
%    =
%    \gamma(x^{\gamma(y)^{-1}}) \gamma(y)
%    =
%    \gamma(y),
%  \end{align*}
%  as $x \in M \le \ker(\gamma)$, and $M$ is characteristic in
%  $(B,\cdot)$. Therefore we have $\gamma(B) = \gamma(Q)$.
%\end{rem}

\begin{prop}\label{prop:duality}
  Let $(G,\cdot)$ be a finite non-abelian group.
  Let $M \ne \{1\}$ be a characteristic subgroup of~$G$ of prime order $p$ such that $M \cap Z(G) = \{1\}$. Let $\lambda\colon a\in G\mapsto \lambda_a\in\operatorname{Aut}(G)$  be a function satisfying $\lambda_{a\cdot\lambda_a(b)}=\lambda_a\lambda_b$, and suppose
  that for every 
  $m \in M$  we have  $\lambda(m) = \iota(m^{-\sigma})$, for  some function
  $\sigma$ on $M$. Then $\sigma \in \operatorname{End}(M) = \mathbb{Z} / p \mathbb{Z}$, and we have that either $\sigma=0$ or $\sigma=1$.
  %\begin{enumerate}
  %\item[\textnormal{(1)}]
  %  either $\sigma = 0$, that is, $M \le \operatorname{Ker}(\lambda)$,
  %\item[\textnormal{(2)}]
  %  or $\sigma =  1$, that is, $\lambda(m) = \iota(m^{-1})$  for $m \in
  %  M$, so that $M \le \operatorname{Ker}(\overline{\lambda})$.
  %\end{enumerate}
\end{prop}

%-m+c+m
%-m+mc
%-m +_o mc=mc-m=m-m+mc-m

\begin{theo}\label{thm:supersylow}
Let $\pi$ be a set of primes and $(B,+,\circ)$ a finite supersoluble skew brace. If $Q$ is any left ideal of $B$, then $Q$ is contained in a Hall $\pi$-sub-skew brace of $B$.
\end{theo}
\begin{proof}
Let $B$ be a counterexample of minimal order, and $M$  a minimal ideal of $B$. By minimality,~$|B/M|$ is a $\pi$-number, and $M$ has prime order $q\not\in\pi$. In particular, $(B,+) = (P,+)+(M,+)$, where~$P$ is a~Hall~\hbox{$\pi$-sub}\-group of
$(B,+)$. If $P$ is normal in $(B,+)$, then $P$ is  characteristic in $(B,+)$, and so it is a Hall $\pi$-sub-skew brace of $B$, a contradiction. Thus, $P$ is not additively normalized by $M$.

For $m\in M$, write \begin{align*}
  \lambda_m
  =
  \begin{bmatrix}
    a_m & b_m\\
    0 & d_m
  \end{bmatrix},
\end{align*}
where $a_m, b_m, d_m$ are defined as in Theorem \ref{thm:semidirect}. Since $M$ is additively characteristic, it is also a subgroup of $(B,\circ)$, so $\lambda(M)$ is a finite~\hbox{$q$-group,} and each $\lambda(m)$
is a $q$-element, for $m \in M$. In particular $a_m$ is an automorphism of order a power of $q$ of the group $M$ of order $q$, and thus $a_m = 1$ for all $m \in M$. We also have $d_m = 1$ for $m \in M$, because $M$ is an ideal.

%as the fact that $M$ is an ideal entails $[\gamma(M), B] \le M$.

In order to obtain the final contradiction, we prove that there exists a function $\tau$ on $M$ such that $b_m(x) = [m^{\tau},x]$, for $x \in P$. The existence of such a function obviously implies that
$\lambda_m = \iota(m^{\tau})$, for $m \in M$, where $\iota$ denotes additive conjugation (on the left). Writing $\tau = -\sigma$, Proposition~\ref{prop:duality} now applies to $G = (B,+)$, so either $\sigma=0$ and $M\leq\operatorname{Ker}(\lambda)$, or $\sigma=1$ and $M$ is easily seen to lay in $\operatorname{Ker}(\lambda^{\operatorname{opp}})$.  In any case, we may apply Lemma~\ref{lemma:ifkerisone} and obtain a contradiction.

\medskip

It remains to prove the existence of such a function $\tau$. We provide two proofs of this fact, one using cohomology results, and the other cohomology-free. As for the former, \hbox{$b_m:P\rightarrow M$} is a derivation for the  action of $P$ on  $M$ by conjugation in $(B,+)$.  Since $P$ and $M$  have coprime orders, we  have $H^{1}(P, M)  = 0$ (see~\cite[III.10.2]{Brown}), that is, all such derivations are inner, so~\hbox{$b_m(x)=[a,x]_+=a+x-a-x$} for some $a\in M$ and for all $x\in P$. This obviously proves our claim.

\medskip

The latter proof is as follows. Let $c_{1}, c_{2} \in C=C_{(P,+)}(M,+)$ and $m\in M$. Then 
\begin{align*}
  b_m(c_{1} + c_{2})
  =
  b_m(c_{1})+ {}^{+,c_1}b_m(c_{2})
  =
b_m(c_{1})+ b_m(c_{2}).
\end{align*}
Therefore ${b_m}_{\restriction C} : (C,+) \longrightarrow (M,+)$ is a group homomorphism, and since $\gcd(|P|, |M|) = 1$, it is trivial. Let $t \in P$ be such that $(P / C,+) =  \langle t +C\rangle_+$. Then
every element of $P$ is of the form $\ell t+c$, for some non-negative integer $\ell$ and some
$c \in C$, so that
\begin{align*}
  b_m(\ell t+c)
  =
  b_m(\ell t)
  =
  {}^{1+t+\ldots +t^{\ell-1}}(b_m(t)).
\end{align*}
%b(t+t)=b(t)+^tb(t)
%b(t+2t)=...
It follows that $b_m$ is completely determined by its value $b_m(t)\in M$, so there are at most $q$ possibilities for $b_m$. But for $m \in M$, the map $\beta_m: x\in P\mapsto [m,x]\in M$ satisfies
\begin{align*}
  \beta_m(x+y)
  =
  [m,x +y]_+
  =
  [m, x]_+ + {}^{x}[m, y]_+
  =
  \beta_m(x)+{}^x\beta_m(y),
\end{align*}
and they are $q$ distinct maps, as for $m, n \in M$ and $x \in P$ we have
\begin{align*}
  [m - n,x]_+
  =
  [m, x]_+ - [n,x]_+,
\end{align*}
as $M$ is abelian, and $[m-n,x]_+ = 0$ for all $x \in P$ if
and only if $m-n=0$, as $P$ acts non-trivially on $M$. This again proves the existence of the map $\tau$.
\end{proof}

\section{An insoluble example}\label{insoluble}

In this section we show that two relevant examples of a simple skew brace, given by~Byott in~\cite{byott}, Theorem 1.1, have Hall sub-skew braces for every set of primes dividing its order. Since each one of these skew braces is the opposite of the other, we may just work with one of them, call it $(B,+,\circ)$. By construction, $(B,+)\simeq P\rtimes C$, where $P$ is an elementary abelian $p$-group, $C$ is a cyclic group of order $q$, and $p,q$ are distinct primes such that $q$ divides $(p^p-1)/(p-1)$; moreover, $(B,\circ)\simeq D\rtimes R$, where $D$ is cyclic of order $q$, and $R$ is a $p$-group. Clearly, $P$ is a~Sy\-low~\hbox{$p$-sub-skew} brace of $G$, and every $p$-sub-skew brace of $B$ is contained in $P$. We need to show that $B$ has sub-skew braces of order $q$. To see this, we consider the natural semidirect product $G=(B,+)\rtimes_\lambda(B,\circ)$. Then we extend the unique~Sy\-low~\hbox{$q$-sub}\-group $D$ of $(B,\circ)$ to a~Sy\-low~\hbox{$q$-sub}\-group $T=S\rtimes_\lambda D$ of $G$. But $T$ has order $p^2$, so is abelian, and hence $D\ast S=\{0\}$. Let $g$ be any element of $(B,\circ)$. Since $D\trianglelefteq (B,\circ)$, we have that $D\ast \lambda_g(S)=\{0\}$. Thus, if we call~$N$ the subgroup of $(B,+)$ generated by the $\lambda_g(S)$, with $g\in (B,\circ)$, then $D\ast N=\{0\}$. But~$N$ is an additive subgroup of $(B,+)$ containing $S$, so either $S$ is $\lambda$-invariant or $N$ contains $B$ because the action of $C$ on $P$ is irreducible. In the former case $S$ is a Sylow $q$-sub-skew brace of $B$, while in the latter, $D$ is a Sylow $q$-sub-skew brace because $D\ast D=\{0\}$. Thus, in any case, $B$ has~Sy\-low~\hbox{$q$-sub-skew} brace.

\bigskip
\bigskip

\noindent {\bf Acknowledgments}\quad  
Ferrara and Trombetti are members of the non-profit association ‘‘AGTA --- Advances in Group Theory and Applications’’ (www.advgrouptheory.com). Funded by the European Union - Next Generation EU, Missione 4 Componente 1 CUP B53D23009410006, PRIN 2022 --- 2022PSTWLB --- Group Theory and Applications. 
Del Corso gratefully acknowledges the support of the MIUR Excellence Department Project awarded to the Department of Mathematics, University of Pisa, CUP I57G22000700001.
Del Corso has performed this activity in the framework of the PRIN 2022, title ‘‘Semiabelian varieties, Galois representations and related Diophantine problems’’. All the authors are supported by \hbox{GNSAGA} (INdAM).

%%%%%%%%%%%%%%%%%%%%

\newpage

\begin{flushleft}
\rule{8cm}{0.4pt}\\
\end{flushleft}

\bigskip
\bigskip

{
\sloppy
\noindent
Andrea Caranti

\noindent
Dipartimento di Matematica

\noindent
Università degli Studi di
Trento

\noindent
via Sommarive 14, I-38123 Trento (Italy)

\noindent
e-mail: andrea.caranti@unitn.it

\noindent
url: https://caranti.maths.unitn.it/
}

\bigskip
\bigskip

{
\sloppy
\noindent
Ilaria Del Corso

\noindent
Dipartimento di Matematica

\noindent
Università di Pisa

\noindent
Largo Bruno Pontecorvo, 5, 56127 Pisa (Italy)

\noindent
e-mail: ilaria.delcorso@unipi.it

\noindent
url: http://people.dm.unipi.it/delcorso/
}

\bigskip
\bigskip

{
\sloppy
\noindent
Massimiliano di Matteo

\noindent
Dipartimento di Matematica e Fisica

\noindent
Università degli Studi della Campania  ``Luigi Vanvitelli''

\noindent
viale Lincoln 5, Caserta (Italy)

\noindent
e-mail: massimiliano.dimatteo@unicampania.it 
}

\bigskip
\bigskip

{
\sloppy
\noindent
Maria Ferrara

\noindent
Dipartimento di Ingegneria

\noindent
Facoltà di Ingegneria e Informatica

\noindent
Università Pegaso

\noindent
e-mail: maria.ferrara1@unipegaso.it 

}

\bigskip
\bigskip

{
\sloppy
\noindent
Marco Trombetti

\noindent 
Dipartimento di Matematica e Applicazioni ``Renato Caccioppoli''

\noindent
Università di Napoli Federico II

\noindent
Complesso Universitario Monte S. Angelo

\noindent
Via Cintia, Napoli (Italy)

\noindent
e-mail: marco.trombetti@unina.it 

}

%\bibliographystyle{plain}
%\bibliography{sample}

\end{document}